\documentclass[bj,authoryear]{imsart}

\usepackage{wrapfig}
\usepackage{tensor}
\usepackage{caption}

\startlocaldefs
\theoremstyle{plain}
\newtheorem{theo}{Theorem}
\newtheorem{propo}[theo]{Proposition}

\newtheorem{lemma}[theo]{Lemma}

\theoremstyle{remark}
\newtheorem{remark}[theo]{Remark}

\newcommand{\eps}{\varepsilon} 
\newcommand{\DD}{\mathcal D}
\newcommand{\E}{\mathbb{E}}

\newcommand{\Cov}{\mathrm{Cov}}
\newcommand{\Corr}{\mathrm{Corr}}

\newcommand{\Var}{\mathbb{V}\mathrm{ar}}

\newcommand{\wh}{\widehat}

\newcommand{\oc}[1]{\overset{\circ}{#1}} 

\newcommand{\bal}[1]{\begin{align*}#1\end{align*}}
\newcommand{\beq}[1]{\begin{equation}#1\end{equation}}

\endlocaldefs

\begin{document}

\begin{frontmatter}
\title{On the joint distribution of the area and the number of peaks for Bernoulli excursions}
\runtitle{Joint distribution of the area and the number of peaks}

\begin{aug}
\author[A]{\inits{V}\fnms{Vladislav}~\snm{Kargin}\ead[label=e1]{vkargin@binghamton.edu}}
\address[A]{Department of Mathematics and Statistics, Binghamton University, USA\printead[presep={,\ }]{e1}}

\end{aug}

\begin{abstract}
Let $P_n$ be a random Bernoulli excursion of length $2n$. We show that the area under $P_n$ and the number of peaks of $P_n$ are asymptotically independent. We also show that these statistics have the correlation coefficient asymptotic to $c /\sqrt{n}$ for large $n$, where $c < 0$, and explicitly compute the coefficient $c$.
\end{abstract}


\end{frontmatter}


\section{Introduction} 

A Bernoulli walk of length $2n$ is a sequence of integers $\eta_i$, $i = 0, \ldots, 2n$, such that $\eta_0 = 0$ and $\eta_i - \eta_{i - 1} \in \{-1, 1\}$ for $i \geq 1$. A Bernoulli excursion has an additional requirement that $\eta_i \geq 0$ for all $i$ and $\eta_{2n} = 0$. A random Bernoulli excursion is chosen uniformly at random from the set of all Bernoulli excursions.

\begin{figure}[htbp]
\centering
 \includegraphics[width=.4\textwidth]{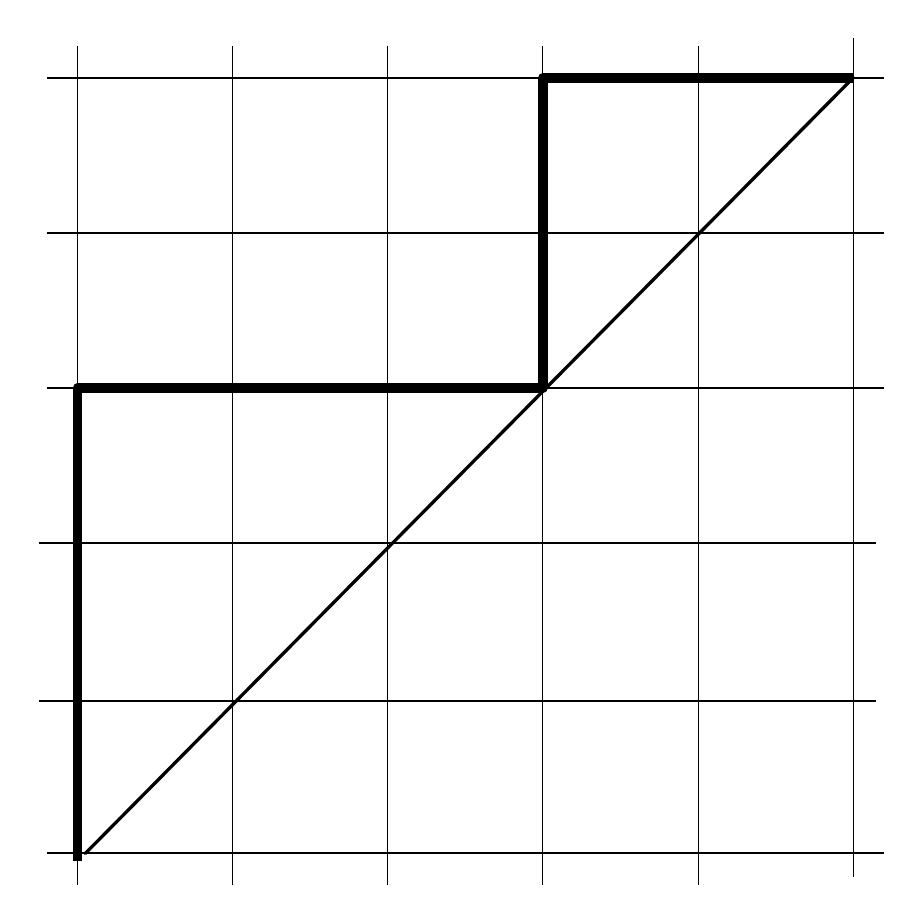}
 \captionsetup{justification=centering}
 \vspace{-1em}
 \caption{{\footnotesize A Dyck path NNNEEENNEE of half-length $n = 5$ with two peaks and area 4.}}
 \label{figDyckPath}
\end{figure}
It will be convenient for us to use the alternative language of Dyck lattice paths, which correspond bijectively to  Bernoulli excursions. 
Recall that a \emph{Dyck path} of half-length $n$ is a lattice path from $(0, 0)$ to $(n, n)$ consisting of $n$ horizontal steps ``East'' from $(i, j)$ to $(i + 1, j)$ and $n$ vertical steps ``North'' from $(i, j)$ to $(i, j + 1)$, such that all points on the path satisfy $i \leq j$, that is, the path lies on or above the line $y = x$. See an example in Figure \ref{figDyckPath}. The bijection with Bernoulli excursions is given by the map: ``North'' and ``East'' at step $i$ correspond to $ \eta_i - \eta_{i - 1} = +1, -1$, respectively.

We will denote the set of all Dyck paths of half-length $n$ by $\DD_n$.  The size of this set is given by the Catalan number 
\bal{
C_n = \frac{1}{n + 1} \binom{2n}{n}.
}

A \emph{peak} in a Dyck path $P$ corresponds to a subsequence  $NE$. The \emph{roughness} of path $P$ is defined as the number of peaks of $P$, which we denote  $pk(P)$.  Another quantity of interest is the \emph{area} of a Dyck path $P$, $a(P)$, which we define as the number of whole unit squares below the path $P$ and above the line $y = x$. For the path in Figure \ref{figDyckPath}, $pk(P) = 2$ and $a(P) = 4$.

Now suppose that the Dyck paths are sampled uniformly at random from $\DD_n$ and define $X_n$ and $Y_n$ as the area and the number of peaks for a random Dyck path, respectively.

An explicit formula for the distribution of $Y_n$ was found by \cite{narayana59} and the asymptotic normality of  $Y_n$ follows from the results of \cite{kolchin84} (Theorem 2.3.1) for random graphs; the asymptotic normality and convergence of all moments was proven by a different method in \cite{janson2001} (Example 3.4). Yet another method of proof can be found in \cite{drmota2009} (Theorem 3.13). Our Proposition \ref{propoAlpha} gives another proof of this convergence by an explicit calculation of the moment limits. So, for large $n$, the distribution of $Y_n$ is approximately normal with expectation $\sim n/2$ and variance $\sim n/8$.

For the random area $X_n$,   Theorem 5 in \cite{takacs1992} implies that $X_n/n^{3/2}$ has a limiting distribution equal to the distribution of $\sqrt{2} \mathcal{B}_{ex}$, where  
\bal{
\mathcal{B}_{ex} = \int_0^1 B_{ex}(t) \, dt, 
}
and $B_{ex}(t)$ is a normalized Brownian excursion. 
The distribution of $\mathcal{B}_{ex}$ was studied in \cite{louchard1984} and \cite{louchard1984b}, where it was described in terms of the Airy function. This distribution is often called the Airy distribution or the Airy distribution of the area type (see, for example, Section 2 in \cite{janson2007} and  Proposition VII.15 on p. 534 in \cite{flajolet_sedgewick2009}). See also \cite{takacs91}, \cite{flajolet_louchard2001} for further analytical properties of the Airy distribution and a survey paper of \cite{janson2007} for additional details about results in this area of research.

Recently, there has been a growing interest  in the joint distribution of random statistics of combinatorial structures. For example, \cite{the2004} investigated the joint distribution of the area and the inertial moment for Bernoulli  bridges and excursions, with the inertial moment defined as $\sum_{i = 0}^n x_i^2.$  \cite{richard2009} generalized \citeauthor{the2004}'s results to the case of $M + 1\geq 2$ moments $u_k = \sum_{i = 0}^n x_i^k$ for $k = 0, \ldots , M$.  \cite{blanco_petersen2014}  studied the joint distribution of the area under a random Dyck path and the rank of the corresponding partition in the lattice of non-crossing partitions.  \cite{chassaing_marckert_yor2000} and \cite{janson2008} studied the joint distribution of the height and width of a random planar tree.

We focus on the \emph{joint distribution} of the random variables $X_n$ and $Y_n$, when $n$ is large. As our first result, we find that $X_n$ and $Y_n$ are asymptotically uncorrelated. More precisely, their correlation coefficient is $-c n^{-1/2} + O(n^{-1})$, where 
\bal{
c = \Big(\frac{20}{3\pi} - 2\Big)^{-1/2} \approx 2.8622\ldots
} 

  If the rescaled limit of $X_n$ and $Y_n$ were a joint Gaussian distribution, then the asymptotically vanishing correlation coefficient between $X_n$ and $Y_n$ would imply the asymptotic independence, that is, the rescaled limit would be a product of marginal Gaussian distributions. In our case, however, the rescaled limit of $X_n$ is the Airy distribution, so in order to investigate the rescaled limit of the joint distribution of $X_n$ and $Y_n$, we check the behavior of correlations between $(X_n)^k$ and $(Y_n)^l$, for all $k, l \geq 1$.  
We prove that these correlations are indeed asymptotically vanishing for all $k, l \geq 1$ and in this way we establish the asymptotic independence of $X_n$ and $Y_n$.

These results can be compared with findings of \cite{labarbe_marckert2007}, who considered a random Bernoulli excursion of length $2n$ conditioned to have a specific fraction of peaks. They assume that  $Y_n/n \to \alpha \in (0, 1)$ as $n \to \infty$. Then, it can be inferred from their results that the excursions with larger $\alpha$ have smaller area. 

In other words,  the results of \citeauthor{labarbe_marckert2007} imply that $\E(n^{-\frac{3}{2}} X_n | Y_n = \lfloor\alpha n \rfloor)$ is decreasing in $\alpha$ for large $n$, while we find that $X_n$ and $Y_n$ are asymptotically independent after rescaling. 

We can suggest the following intuitive explanation for this difference in results. The theorems of \cite{labarbe_marckert2007} imply that a change in $Y_n$ on the scale of $n$ is associated to a change in $X_n$ on the scale of $n^{3/2}$, which is the natural order of fluctuations for $X_n$. However, since $Y_n$ is concentrated, its natural order of fluctuations is $n^{1/2}$. If we assume that the dependence between $X_n$ and $Y_n$ found in \cite{labarbe_marckert2007} operates also on smaller scales, then the fluctuations in $Y_n$ on the order of $n^{1/2}$  are likely to be associated with changes in $X_n$ on order of $n$ in the conditional expectation, and this is smaller than $n^{3/2}$, the natural order of fluctuations in $X_n$. Thus, these fluctuations will not be strong enough to prevent the asymptotic independence after rescaling.

For the proof of our results, we use the method of generating function similar to the method employed by \cite{the2004}. However, in our case, the sub-leading terms of the asymptotic expansion of $\E(X_n^k Y_n^l)$ are also needed due to the fact that random variables $X_n$ and $Y_n$ are of different type: while $Y_n$ concentrates with the growth in $n$,  $X_n$ does not. In order to overcome this difficulty, we will use the ``dominant balance'' method for $q$-functional equations, which was developed by Richard and collaborators in a series of papers: \cite{richard_guttmann2001}, \cite{richard_guttmann_jensen2001}, \cite{richard2002}, \cite{richard_jensen_guttmann2008}, \cite{srt2010}.  The original goal of these papers was to study the lattice polygon models and to find an approach to the enumeration problem for closed self-avoiding walks; however, this method is quite useful in other combinatorial problems, and we employ it to prove the asymptotic independence of $X_n$ and $Y_n$.
 
Finally, we would like to mention that the random variables $X_n$ and $Y_n$ appear also in other combinatorial problems.  Indeed, Bernoulli excursions and Dyck paths are examples of \emph{Catalan structures}, -- the families of combinatorial objects which are enumerated by Catalan numbers. These structures are inter-connected by numerous bijections. A list of these structures is given in Volume 2 of \cite{stanley99}. See also \cite{stanley2015} and the ``Catalan Addendum'' on Stanley's webpage \url{https://math.mit.edu/~rstan/ec/} for a larger list containing 207 Catalan structures. 

For example, one of the bijections between Catalan structures relates Dyck lattice paths of half-length $n$ and planar rooted trees on $n + 1$ vertices. This bijection is given by the depth-first search (``DFS'') walk on the tree.  Under this bijection,  $X_n$ = ``total path length of the tree - $n$'', and $Y_n$ corresponds to the number of leaves in the tree. By using this bijection, Tak\'acs gave the formula for the distribution of the total path length for a random tree (see \cite{takacs94}).  Our results can be similarly reformulated as results about the joint  distribution of the number of leaves and the total path length of a random planar tree.

 In the rest of the paper, we give precise statements of our results in Section \ref{result_statements} and describe the general outline of  proofs in Section \ref{section_strategy}.  The proofs of theorems will be presented in the following order. First, we prove Theorem \ref{mainResult} in Section \ref{proofOfMainResult}. Then we prove some additional results necessary for the proof of Theorem \ref{theoMainResult0} in Section \ref{preliminary_results} and prove a particular case of Theorem \ref{theoMainResult0} in Section \ref{proofOfMainResult2}. Finally, we prove Theorem \ref{theoMainResult0}  in full generality in Section \ref{proofOfAsymptoticIndependence}.


\section{Results}
\subsection{Statements}
\label{result_statements}
Let $\Corr(\xi, \eta)$ denote the correlation coefficient between random variables $\xi$ and $\eta$,  
\bal{
\Corr(\xi, \eta) := \frac{\Cov(\xi, \eta)}{\sqrt{\Var(\xi)\Var(\eta)}}.
}

\begin{theo}
\label{mainResult} Let $X_n$ and $Y_n$ be the area and the roughness of a uniformly random Dyck path $P \in \DD_n$. Then, as $n \to \infty$,
\bal{
\Corr(X_n, Y_n) = - \frac{1}{\sqrt{2}\sqrt{\frac{10}{3\pi} - 1}} n^{-1/2} + O\big(n^{-1}\big).
}
\end{theo}

Numerically, the coefficient is $-2.8622$. As an immediate consequence, we find that $X_n$ and $Y_n$ are asymptotically uncorrelated for large $n$. 
 We will prove this theorem in Section \ref{proofOfMainResult}.

Theorem \ref{mainResult} suggests that $X_n$ and $Y_n$ might be asymptotically independent, that is, that the joint distribution of appropriately scaled $X_n$ and $Y_n$ converges to the product of the Gaussian and Airy distributions. 

To investigate this possibility, we note that in order to show the asymptotic independence, it is enough to show that as $n \to \infty$,
\beq{
\label{asympIndependence}
\E\Big[g\Big(\frac{X_n - \E X_n}{n^{3/2}}\Big)h\Big(\frac{Y_n - \E Y_n}{n^{1/2}}\Big)\Big] - \E g\Big(\frac{X_n - \E X_n}{n^{3/2}}\Big) \E h\Big(\frac{Y_n - \E Y_n}{n^{1/2}}\Big) \to 0,
} 
 for all bounded continuous functions $g$ and $h$. It is known that the scaling limit distributions for $X_n$ and $Y_n$ have exponentially-declining tails. (For $Y_n$, this is the Gaussian distribution, and for $X_n$, the scaling limit is the Airy distribution, which also has thin tails -- see the proof of Theorem 3.1. in \cite{csorgo_shi_yor99} and Example 4.2 in \cite{fill_janson2009}.) Since the polynomial functions are dense among the continuous functions on compact supports,  we can use truncation to show that it is enough to check (\ref{asympIndependence}) for polynomial $g$ and $h$.

In particular, the asymptotic independence will follow, if we show that the correlations between $(X_n)^k$ and $(Y_n)^l$ converge to zero for all $k, l \geq 1$, as $n \to \infty$. Let 
\begin{equation}
\label{correlations}
  \Corr(X_n^k, Y_n^l)  = \frac{\E(X_n^k Y_n^l) - \E(X_n^k) \E(Y_n^l)}
 {\Var(X_n^k)^{1/2} \Var(Y_n^l)^{1/2}}
\end{equation}
 and 

 \begin{equation}
 \label{normalization}
 \wh X_n = \frac{X_n}{\sqrt{2} n^{3/2}} \text{ and }  \wh Y_n = \frac{Y_n - \E Y_n}{\sqrt{\Var(Y_n)}}.
\end{equation}
 
   \begin{theo} 
 \label{theoMainResult0}
 For all integer $k, l \geq 1$, as $n \to \infty$, 
 \bal{
 \Corr(X_n^k, Y_n^l)  = O\big(n^{-1/2}\big),
 }
 and $(\wh X_n, \wh Y_n)$ converges in distribution to the product of the Airy distribution (i.e., the distribution of $\mathcal{B}_{ex}$) and the standard normal distribution. 
 \end{theo}
\begin{remark}
While the normalization of $Y_n$ by $\E Y_n$ in (\ref{normalization}) is essential, the normalization of $X_n$ by $\E X_n$ can be omitted since $\E X_n$ and $\sqrt{\Var(X_n)}$ have the same order $n^{3/2}$ as $n \to \infty$. 
\end{remark}

 \subsection{Strategy of the proof}
 \label{section_strategy}
 
 The statistics $X_n$ and $Y_n$ are in a certain sense very different from each other. Namely, as $n$ grows, $Y_n$ (the roughness of the path) concentrates around its mean, that is, its  standard deviation becomes negligible relative to its mean: $\E(Y_n) \asymp n$ (that is, $\E(Y_n) \sim c n$ for a positive constant $c$), and $\sqrt{\Var(Y_n)} \asymp n^{1/2}$. In contrast, $X_n$ does not concentrate: $\E(X_n) \asymp n^{3/2}$ and $\sqrt{\Var(X_n)} \asymp n^{3/2}$ as well. 
 
 This observation can be extended to higher moments: $\E(Y_n^l) \asymp n^l$, and $\sqrt{\Var(Y_n^l)} \asymp n^{l - 1/2}$, while $\E(X_n^k)  \asymp \sqrt{\Var(X_n^k)}  \asymp n^{3k/2}$.

  This implies that the denominator in (\ref{correlations}) is $n^{s}$ where $s =\frac{3}{2}k + l - \frac{1}{2}$, and we aim to show that the the numerator is $O\big(n^{s - \frac{1}{2}}\big)$. 
  
Our approach to this problem will be to write the expansions:
\begin{align}
\notag
\E( X_n^k) &\sim a_k^{(0)} n^{3k/2} + a_k^{(1)}n^{(3k - 1)/2} + \ldots, 
\\
\notag
\E( Y_n^l) &\sim b_l^{(0)} n^{l} + b_l^{(1)}n^{l - 1/2} + \ldots,
\\
\label{expansion_moments}
\E( X_n^k Y_n^l) &\sim c_{k,l}^{(0)} n^{3k/2 + l} + c_{k,l}^{(1)}n^{(3k - 1)/2 + l} + \ldots .
\end{align}  
(Superscript $(0)$ is for the leading order and $(1)$ is for the first subleading order.) Then, we have
\bal{
\E( X_n^k Y_n^l) - \E( X_n^k)\E( Y_n^l) &\sim [c_{k,l}^{(0)} - a_k^{(0)} b_l^{(0)}] n^{3k/2 + l} 
\\
&+ [c_{k,l}^{(1)} - a_k^{(0)}  b_l^{(1)} -  a_k^{(1)} b_l^{(0)}] n^{(3k - 1)/2 + l} + \ldots
}  
From this asymptotic expansion, we can conclude that it must be that $c_{k,l}^{(0)} - a_k^{(0)} b_l^{(0)} = 0$. (Otherwise the absolute value of the correlation $|Corr(X_n^k, Y_n^l)|$ would exceed 1 for sufficiently large $n$.) Hence, we only need to show that $c_{k,l}^{(1)} - a_k^{(0)}  b_l^{(1)} -  a_k^{(1)} b_l^{(0)} = 0$. 

One step in the proof is to show that $b_l^{(1)} = 0$ for all $l \geq 1$. 
(See Remark \ref{remark_justification} after Proposition  \ref{propoAlpha} for the justification of this identity.) After this identity is proved, in order to show that $Corr(X_n^k, Y_n^l) = O(n^{-1/2})$, we need to show that 
\begin{equation}
\label{main_identity}
c_{k,l}^{(1)} = a_k^{(1)} b_l^{(0)}.
\end{equation}

The coefficients in expansions (\ref{expansion_moments}) can be obtained by relating these coefficients (``moments'') to the asymptotic expansion of the moment generating function in a neighborhood of the smallest singularity. This step is followed by the method of dominant balance, which is a powerful technique to perform “moment pumping” systematically and calculate the coefficients in the asymptotic expansion of  the moment generating function (``mgf''). This is achieved by an appropriate rescaling of the mgf variables, so that the functional equation satisfied by the mgf is converted into differential equations for the generating functions of the leading and sub-leading singular parts of the mgf.

In more detail, we define
\beq{
\label{def_D}
D(q, u, z) = \sum_{P \in \DD} q^{a(P)}(1 + u)^{pk(P)} z^{|P|},
}
where the sum is over all Dyck paths, $|P|$ is the half-length of path $P$, and where by convention we set $D(q, u, 0) = 1$. 

The joint factorial moments of $X_n$ and $Y_n$ can be obtained from the partial derivatives of $D(q, u, z)$ with respect to $q$ and $u$, evaluated at the point $(q, u) = (1,0)$,

\begin{align}
\label{moments_derivatives}
&\E\Big[X_n (X_n - 1) \ldots (X_n - k + 1) Y_n (Y_n - 1) \ldots (Y_n - l + 1)\Big] \notag
\\
&
= \frac{1}{C_n}\Big\{[z^n]\frac{\partial^k}{\partial q^k} \frac{\partial^l}{\partial u^l} D(q, u, z)\Big|_{q = 1, u = 0} \Big\},
\end{align} 
where we use notation $[z^n] f(z)$ for the coefficient of the $n$-th power of $z$ in the series for $f(z)$, and where $C_n$ is the Catalan number (i.e., the total number of Dyck paths on $2n$ steps).

It is possible to show by induction that the first two highest terms in the asymptotic expansion of $\E(X_n^k Y_n^l)$ coincide with the corresponding terms in the expansion of the factorial moment, so we concentrate on evaluation of $[z^n]\frac{\partial^k}{\partial q^k} \frac{\partial^l}{\partial u^l} D(q, u, z)\Big|_{q = 1, u = 0}$. In order to get coefficients in (\ref{expansion_moments}), we need the asymptotics of these expressions for large $n$.

The \emph{transfer method} is a way to get the asymptotics of $[z^n] f(z)$ for a given $f(z)$. The first step of this method asks for an asymptotic expansion of $f(z)$ in powers of $(1 - z/z_0)$, where $z_0$ is the singularity of $f(z)$ with the smallest absolute value. (In our case, $z_0 = 1/4$.) Then, one can translate the coefficients in this asymptotic expansion to asymptotics of $[z^n] f(z)$ by using a set of specific transfer formulas. See  Section VI in \cite{flajolet_sedgewick2009} for a detailed description of this method.

In particular,  near $z = \frac{1}{4}$, we have the asymptotic expansion 
\begin{align}
\label{exp_derivative}
\frac{\partial^k}{\partial q^k} \frac{\partial^l}{\partial u^l} D(q, u, z)\Big|_{q = 1, u = 0}  &= \frac{\alpha_{k,l}}{(1 - 4z)^{3k/2 + l - 1/2}} - \frac{\beta_{k,l}}{(1 - 4z)^{3k/2 + l - 1}}\\
& + O\big((1 - 4z)^{-(3k/2 + l - 3/2)} \big),  \notag
\end{align}
for $k, l \geq 0$. 
This expansion can be derived by using induction.  The induction starts with the cases $(k, l) \in \{(0, 0), (1, 0), (0, 1)\}$, which we consider in Section \ref{proofOfMainResult}, and then proceeds by the route taken in the proof of expansion in Lemma 3.5 in \cite{srt2010} by differentiating the functional equation for $D(q, u, z)$. In particular, the induction shows that the order of the singularity of the function $\frac{\partial^k}{\partial q^k} \frac{\partial^l}{\partial u^l} D(q, u, z)\Big|_{q = 1, u = 0}$ at $z = \frac{1}{4}$ is $\frac{3}{2}k + l - \frac{1}{2}$. This could be also guessed from the known results on the asymptotics of moments of random variables $X_n$ and $Y_n$ for large $n$.

 Then, by Theorems VI.1 and VI.4 in \cite{flajolet_sedgewick2009}, if either $k > 0$ or $l >1$, we have 
 \begin{align}
 \label{expansion_zn}
[z^n]\frac{\partial^k}{\partial q^k} \frac{\partial^l}{\partial u^l} D(q, u, z)\Big|_{q = 1, u = 0} &\sim 4^n\Big[ \alpha_{k,l}\frac{n^{3k/2 + l - 3/2}}{\Gamma(3k/2 + l - 1/2)}  + \beta_{k,l}\frac{n^{3k/2 + l - 2}}{\Gamma(3k/2 + l - 1)} \\
& + O(n^{3k/2 + l - 5/2})\Big]. \notag
\end{align}
(The case $k = 0$ and $l = 1$ is special, since in this case $3k/2 + l - 1 = 0$ and the second term in the expansion (\ref{exp_derivative}) is not singular. In this case, the second term in the expansion (\ref{expansion_zn}) should be replaced by zero.)  
According to Figure VI.3 in \cite{flajolet_sedgewick2009}, 
\bal{
C_n \sim \frac{4^n}{\sqrt{\pi n^3}} \Big(1 - \frac{9}{8} n^{-1} + O(n^{-2})\Big),
}
and we get from (\ref{expansion_moments}), (\ref{moments_derivatives}) and (\ref{expansion_zn}) that
\bal{
c_{k,l}^{(1)} = \frac{\sqrt \pi \beta_{k, l}}{\Gamma(3k/2 + l - 1)},\
a_{k}^{(1)} = \frac{\sqrt \pi \beta_{k, 0}}{\Gamma(3k/2 - 1)}, \
b_{l}^{(0)} = \frac{\sqrt \pi \alpha_{0, l}}{\Gamma(l - 1/2)}.
}
Hence, in order to prove (\ref{main_identity}), we need to obtain the relation
\beq{
\label{main_identity2}
\beta_{k, l} = \frac{\sqrt \pi \Gamma(3k/2 + l - 1)}{\Gamma(3k/2 - 1) \Gamma(l - 1/2)} \beta_{k, 0} \alpha_{0, l}.
}
In Proposition  \ref{propoAlpha} below we show that for every $l \geq 1$,
\beq{
\alpha_{0, l} = \frac{\Gamma(l - 1/2)}{2^l \sqrt \pi},
}
so (\ref{main_identity2}) reduces to 
\beq{
\label{main_identity2a}
\beta_{k, l} = \frac{\Gamma(3k/2 + l - 1)}{2^l\Gamma(3k/2 - 1)} \beta_{k, 0} \equiv \frac{(3k/2 - 1)_l}{2^l}\beta_{k, 0},
}
where $(x)_l$ is the Pochhammer symbol, $(x)_l := x(x+1) \ldots (x + l -1)$.

This is where the method of dominant balance contributes. The method of dominant balance is a powerful technique to perform ``moment pumping” systematically, also for sub-leading corrections to the asymptotic behavior. This is achieved by an appropriate rescaling of the generating function variables in the neighborhood of singularity. This turns the given functional equation into a sequence of differential equations. For example, we write $\wh F(s, \eps) = D(q, u, z)|_{u = 0}$ where $1 - q = \eps^3$, $1 - 4z = s \eps^2$. (This choice for the change of variables is suggested by the asymptotic behavior of the moments of  random variables $X_n$ and $Y_n$.) Then, the dominant balance method suggests to search the solution of the functional equation on $D(q, u, z)$ in the form 
 \begin{equation}
 \label{exp_dbalance_intro}
\wh F(s, \eps) = 2 + F_0(s) \eps + F_1(s) \eps^2 + \ldots ,
\end{equation}
 where $F_i(s)$ are Laurent series in powers of $s^{1/2}$. (This specific form is suggested by the solution of the functional equation when $q = 1$ and $u = 0$, see Eqn. (\ref{equ_Fz}) below.) More precisely,
 \begin{equation}
 \label{expansions_F0_F1_0}
 F_0(s) = \sum_{k=0}^\infty \frac{f_{k}^{(0)}}{s^{3k/2 - 1/2}} \text{ and } F_1(s) = \sum_{k=0}^\infty \frac{f_{k}^{(1)}}{s^{3k/2 - 1}}, 
\end{equation}
and these expansions are closely related to the asymptotic expansion (\ref{exp_derivative}). In particular, we will show that 
 \bal{
 (-1)^k k ! f_{k}^{(0)} = \alpha_{k,0} \text{ and } (-1)^k k! f_{k}^{(1)} = - \beta_{k,0}.
 }
 
 The expansions (\ref{exp_dbalance_intro}) and  (\ref{expansions_F0_F1_0}) are formal power series, and they are especially suitable for solving the functional equation on $D(q, u, z)$.  In particular, (\ref{expansions_F0_F1_0}) are generating series for coefficients $\alpha_{k, 0}$ and $\beta_{k, 0}$. We follow \cite{srt2010} here. 

By using the expansion (\ref{exp_dbalance_intro}) and the functional equation for $D(q, u, z)$, we can obtain the differential equations for functions $F_0(s)$, $F_1(s)$, \ldots. For example, we have 
\bal{
F'_0(s) &= 2s - \frac{1}{2} \Big(F_0(s)\Big)^2.
}
This equation implies a recursion on coefficients $f_{k}^{(0)}$ which can be used to derive the Tak\'acz recursion for the main asymptotic term of the area distribution moments. 

Our application of the method of dominant balance follows these main ideas, however, it is more involved.  In order to illustrate the main ideas of the method, we first prove Theorem  \ref{theoMainResult0} for the special case $k \geq 1$ and $l = 1$ in Section \ref{proofOfMainResult2}. This involves differential equations for $F_0(s)$ and $F_1(s)$ and differential equations for $G_0(s)$ and $G_1(s)$, which are Laurent series in the dominant balance approximation for $\wh G(s, \eps) = \partial_u D(q, u, z)|_{u = 0} - 1$. From these differential equations we will find that $G_0(s) = - \frac{1}{2} F'_0(s)$ and $G_1(s) = \frac{1}{2} - \frac{1}{2}F'_1(s)$ and show in (\ref{main_identity4bis}) that this implies that
\beq{
\label{main_identity4}
\beta_{k,1} = \frac{3k - 2}{4}\beta_{k,0},
}
which is (\ref{main_identity2a}) for $l = 1$.

In Section \ref{proofOfAsymptoticIndependence}, we consider general $k$ and $l$. Here we use a somewhat more general change of coordinates, $1 - q = \eps^3$, $1 - 4z = s \eps^2$, $u = t \eps^2$ and define $\wh F(s, t, \eps) = D(q, u, z) - u$.  We look for approximation of $\wh F(s, t, \eps)$ in the form 
\bal{
\wh F(s, t, \eps) = 2 + (2s - t/2) \eps^2 + \eps\Big(F^0(s, t)+ F^1(s,t) \eps + F^2(s, t) \eps^2 + \ldots \Big).
}
Together with the functional equation on $D(q, u, z)$, this gives some partial differential equations on $F^0(s,t)$, $F^1(s,t)$, \ldots, which we can use to derive relations between coefficients in the Laurent expansions of functions $F^0(s,t)|_{t = 0}$ and $F^1(s,t)|_{t = 0}$. Eventually, this allows us to prove the identity (\ref{main_identity2a}) and completes the proof of Theorem \ref{theoMainResult0}.


\section{Proofs}

The basis for our calculations is Theorem 6.14 in \cite{petersen2015}, which implies the following formula:
\begin{equation}
\label{continuous fraction}
 D(q, u, z) = u +\cfrac{1}{1 - z u - \cfrac{z}{1 - q z u -  \cfrac{q z}{1 - 
q^2 z u - \cfrac{q^2 z}{\ldots}}}}.
\end{equation}
(The formula is slightly adjusted since we use peaks of Dyck paths and \citeauthor{petersen2015} uses valleys in the definition of the generating function.)

In fact, it will be more convenient to work with a slightly different function
\beq{
\label{eqDandF}
 F(q, u, z) := D(q, u, z) - u.
}
For $F(q, u, z)$, equation (\ref{continuous fraction}) means that we have the following functional equation:
\beq{
\label{eqFun}
F(q, u, z) = \frac{1}{1 - zu - zF(q, u, q z)}.
}

By differentiating, we can also obtain functional equations for derivatives, for example:
\begin{eqnarray}
F_u(q, u, z) &= \frac{z}{\big(1 - zu - zF(q, u, q z)\big)^2}\big(1 + F_u(q, u, q z)\big) \notag
 \\
 \label{eqFu}
& = z [F(q, u, z)]^2\big(1 + F_u(q, u, q z)\big).
\end{eqnarray}
\begin{remark} The function $F(q, 0, z)$ is the generalized Rogers-Ramanujan function (with the reversed sign for the argument $z$). The standard (and famous) Rogers-Ramanujan function is $F(q, 0, -1)$. See \cite{berndt_yee2003} for additional information. 
\end{remark}

\subsection{Proof of Theorem \ref{mainResult}}
\label{proofOfMainResult}
In addition to formula (\ref{eqFu}), we can derive formulas for other derivatives of $F(q, u, z)$. For conciseness, we suppress the first two arguments of $F(q, u, z)$ in the following formulas. For example, $
F(z) := F(q, u, z)$ and $F(qz):= F(q, u, qz)$. We have
\begin{align}
\label{eqFq}
F_q(z) &= z F^2(z) \Big(F_q(qz) + z F_z(qz)\Big),
\\
\label{eqFqq}
F_{qq}(z) &= 2 z F(z) F_q(z) \Big(F_q(qz) + z F_z(qz)\Big)
\\
\notag &+  z F^2(z)  \Big(F_{qq}(qz) + 2z F_{qz}(qz) + z^2 F_{zz}(qz)\Big),
\\
\label{eqFqu}
F_{qu}(z) & = 
2 z F(z) F_u(z)\Big(F_q(qz) + z F_z(qz)\Big) 
\\
\notag &+ z F^2(z)\Big(F_{qu}(qz) + z F_{zu}(qz)\Big),
\\
\label{eqFuu}
F_{uu}(z) &= 2 z F(z) F_u(z) \Big(1 +  F_u(qz)\Big) + z F^2(z) F_{uu}(qz).
\end{align}

Now let us introduce another piece of notation. For any function $f(q, u, z)$, we will write:
\beq{
\label{defi_oc}
\oc f(z) := f(q, u, z)\Big|_{q = 1, u = 0}. 
}
Our generating functions $f(q, u, z)$ have a singularity at $(q, u, z) = (1, 0, 1/4)$ and we will expand $\oc f(z)$ in power series of $(z - 1/4)$.

By using functional equation (\ref{eqFun}), we get 
\beq{
\label{equ_Fz}
\oc F(z) = \frac{1 - \sqrt{1 - 4z}}{2z} = 2\big( 1 - (1 - 4z)^{1/2} + (1 - 4z) - (1 - 4z)^{3/2} + \ldots\big). 
}
and we can calculate $\oc F_z(z)$, $\oc F_{zz}(z)$, $\ldots$ by differentiation. We can also recursively calculate the derivatives of $F(q, u, z)$ with respect to $q$ and $u$, evaluated at $q = 1$ and $u = 0$, by using the equations from the system above. For example, for $\oc F_q(z)$, equation  (\ref{eqFq}) implies that
\beq{
\oc F_q(z) = z \oc F(z)^2\Big(\oc F_q(z) + z \oc F_z(z)\Big),
}
which can be solved as 
\begin{align}
\oc F_q(z) &= \frac{1 - 2z - \sqrt{1 - 4z}}{(1 - 4z) (1  + \sqrt{1 - 4z})} 
\\
\label{expFq_in_z}
& = \frac{1}{2}(1 - 4z)^{-1} - \frac{3}{2}(1 - 4z)^{-1/2} + 2 
- 2 (1 - 4z)^{1/2} + O\big( ( 1 - 4z)\big).  
\end{align} 
The transfer theorem VI.4 from \cite{flajolet_sedgewick2009} allows us to get the asymptotic expansion for the coefficient before $z^n$ in the expansion of $\oc F_q(z)$ in powers of $z$:
\beq{
\label{expFq_in_n}
[z^n] \oc F_q(z) = 4^n \Big[ \frac{1}{2} - \frac{3}{2} \frac{1}{\sqrt{\pi} n^{1/2}} + \Big(\frac{3}{16} + 1\Big) \frac{1}{\sqrt{\pi} n^{3/2}}
+ O\big( n^{-5/2}\big) \Big].
}
In addition, by using Table VI.3 in \cite{flajolet_sedgewick2009}, we write the asymptotic expansion for the Catalan numbers as
\beq{
C_n = \frac{4^n}{\sqrt{\pi n^3}} \Big[ 1 - \frac{9}{8} n^{-1}+ \frac{145}{128} n^{-2} - \frac{1155}{1024} n^{-3} +
+ O\big( n^{-4}\big) \Big].
}
This leads to 
the expansion for the expected value of $X_n$: 
\beq{
\label{EXn}
\E X_n = \frac{[z^n] \oc F_q(z)}{C_n} = \frac{\sqrt{\pi}}{2} n^{3/2} - \frac{3}{2} n + \frac{9 \sqrt{\pi}}{16} n^{1/2} - \frac{1}{2} + O(n^{-1/2}).
}

In a similar fashion, we can calculate the expansion of $\oc F_{qq}(z)$:
\beq{
\label{expFqq_in_z}
\oc F_{qq}(z) = \frac{5}{8}(1 - 4z)^{-5/2} - \frac{3}{2}(1 - 4z)^{-2} 
+ \frac{3}{4} (1 - 4z)^{-3/2} + O\big( ( 1 - 4z)^{-1}\big), 
} 
which implies that the second fractional moment of $X_n$ is 
\bal{
\E \big[X_n (X_n - 1)] = \frac{[z^n] \oc F_{qq} (z)}{C_n}
= \frac{5}{6} n^3  - \frac{3 \sqrt{\pi}}{2} n^{5/2} + 4 n^2 + 
O\big(n^{3/2}), 
}
and since the variance $\Var (X_n) = \E \big[X_n (X_n - 1)] + \E X_n - \big(\E X_n\big)^2$, we  can calculate that
\beq{
\label{VarXn}
\Var (X_n) =  \Big(\frac{5}{6} -\frac{\pi}{4}\Big)n^3  - \Big(\frac{7}{4} -\frac{9 \pi}{16}\Big)n^2  +  
O\big(n).
}
(The expansions (\ref{EXn}) and (\ref{VarXn}) are related to formulas in Section 4 of \cite{takacs1992}. In particular, Theorem 4 in \cite{takacs1992} gives the main asymptotic term for the moments of a random variable $\omega_n$ closely related to our $X_n$.)

For the moments of the random variable $Y_n$ and the cross-moments, we need partial derivatives of $D(q, u, z)$ with respect to $u$; so we note that equation (\ref{eqDandF}) implies that
$\oc D_u(z) = \oc F_u(z) + 1 $, $\oc D_{uu}(z) = \oc F_{uu}(z)$, and 
$\oc D_{qu}(z) = \oc F_{qu}(z)$. 

From equations (\ref{eqFu}), (\ref{eqFqu}), and (\ref{eqFuu}), we obtain
\beq{
\label{expFu_in_z}
\oc F_u(z) = \frac{1}{2}(1 - 4z)^{-1/2} - \frac{1}{2}, 
} 

\beq{
\label{expFqu_in_z}
\oc F_{qu}(z) = \frac{1}{4}(1 - 4z)^{-2} - \frac{3}{8}(1 - 4z)^{-3/2} - \frac{1}{8} (1 - 4z)^{-1}
+ O\big((1 - 4z)^{-1/2}\big),  
} 

\beq{
\label{expFuu_in_z}
\oc F_{uu}(z) = \frac{1}{8}(1 - 4z)^{-3/2} - \frac{1}{4}(1 - 4z)^{-1/2} + \frac{1}{8} (1 - 4z)^{1/2}
+ O\big((1 - 4z)\big).  
} 

This leads to the following results:
\begin{align}
\label{E_Yn}
\E Y_n &=  \frac{[z^n] \oc D_{u} (z)}{C_n}
= \frac{1}{2} (n  + 1)  + O\big(n^{-2}), \\
\Var (Y_n) &= \frac{n}{8} + \frac{1}{16}  + O\big(n^{-1}), \\
\E (X_n Y_n) &= \frac{\sqrt{\pi}}{4} n^{5/2} - \frac{3}{4} n^2 - \frac{13\sqrt{\pi}}{32} n^{3/2} - \frac{3}{4} n + O\big(n^{1/2}\big).
\end{align}
Then, using (\ref{EXn}),
\bal{
\Cov(X_n, Y_n) = \E(X_n Y_n) - \E(X_n) \E(Y_n) 
= -\frac{\sqrt{\pi}}{8} n^{3/2} + \frac{1}{4}n + O\big(n^{1/2}\big)
}
and 
\bal{
\Corr(X_n, Y_n) = \frac{\Cov(X_n, Y_n)}{\sqrt{\Var(X_n) \Var(Y_n)}}
= -\frac{1}{\sqrt{2}\sqrt{\frac{10}{3 \pi} - 1}} n^{-1/2} + O(n^{-1}).
}
This completes the proof of Theorem \ref{mainResult}.

%
%

 \subsection{A preliminary result for the proof of Theorem \ref{theoMainResult0}}
 \label{preliminary_results}
 Let $F(u,z) := F(q, u, z)\big|_{q = 1}$. 
  \begin{propo} 
 \label{propoAlpha}
 For $r \geq 2$, we have that
 \beq{
 \label{formula_lemmaHyper}
  \frac{1}{r!}\frac{\partial^r}{\partial u^r}  F(q, u, z)\big|_{q = 1, u = 0}
 = \frac{1}{r!}\frac{\partial^r}{\partial u^r}  F(u, z)\big|_{u = 0} = \frac{z^r P_{r - 2}(z)}{(1 - 4z)^{r - \frac{1}{2}}}, 
 }
 where $P_n(z)$ is a polynomial given by formula 
 \beq{
 \label{defi_P}
 P_n(z) = \tensor*[_2]{F}{_1}\Big(\frac{1 - n}{2}, -\frac{n}{2}, 2, 4z\Big),
 }
and  $\tensor*[_2]{F}{_1}(a, b, c, z)$ is the Gauss hypergeometric function. Moreover, for $r \geq 1$, we have the following asymptotic expansion around $z = 1/4$:
 \beq{
 \label{formula_lemmaHyper_a}
 \frac{\partial^r}{\partial u^r}  F(q, u, z)\big|_{q = 1, u = 0}= \frac{\alpha_{0,r}}{(1 - 4z)^{(2r - 1)/2}} - \frac{\beta_{0,r}}{(1 - 4z)^{r - 1}} + O\big((1 - 4z)^{-\frac{2r - 3}{2}} \big), 
 }
 where $\alpha_{0,1} = \frac{1}{2}$, $\beta_{0,1} = \frac{1}{2}$,
 and for $r \geq 2$, we have 
 \bal{
 \alpha_{0,r} = \frac{\Gamma\big(r - \frac{1}{2}\big)}{2^r \sqrt{\pi}}=\frac{(2r - 3)!!}{2^{2r - 1}}, \quad
 \beta_{0,r} = 0.
 }
 \end{propo}   
 \begin{remark}
 \label{remark_justification} Since $b_l^{(1)}$ in (\ref{expansion_moments}) are proportional to $\beta_{0,l}$ for 
   $l \geq 2$, one implication of Proposition \ref{propoAlpha} is that $b_l^{(1)} = 0$ for $l \geq 2$. One can also check that $b_1^{(1)} = 0$  (see equation (\ref{E_Yn})). This justifies the claim that we made in Section  \ref{section_strategy}.  
 \end{remark}
   
 \begin{remark} Since (i) $\alpha_{0,r}$ and $\beta_{0,r}$ determine the asymptotic behavior of the moments $\E(Y_n^r)$,  and (ii) we can identify $\alpha_{0, r}$ as the moments of a normal distribution, Proposition \ref{propoAlpha} provides us with an alternative proof for the rescaled convergence of $Y_n$ to the standard normal distribution. 
 \end{remark}
   
\begin{proof}[Proof of Proposition \ref{propoAlpha}] Note that since we do not differentiate over $q$, we can set $q = 1$ immediately. From (\ref{eqFun}),
 \begin{equation}
 \label{solution_q_1}
 F(u, z) = \frac{1 - uz - \sqrt{(1 - uz)^2 - 4z}}{2z}.
 \end{equation}
 (The case with the positive sign before the square root is ruled out because the definition of $D(q, u, z)$ in (\ref{def_D}) implies that $\lim_{z \to 0} F(0,z) = 1$.) In order to prove (\ref{formula_lemmaHyper}), write 
\bal{
 F_{uu}(u, z) &= \frac{2 z^2}{ \Big((1 - uz)^2 - 4z\Big)^{3/2}} = \frac{2 z^2}{ \Big((1 - t)^2 - 4z\Big)^{3/2}},
 }
where we changed variables, $t = uz$.  By repeatedly using the identity
\bal{
\frac{\partial f(t,z)}{\partial u} = \frac{\partial f(t,z)}{\partial t} \frac{\partial t}{\partial u} + \frac{\partial f(t,z)}{\partial z} \frac{\partial z}{\partial u} = z \frac{\partial f(t,z)}{\partial t},
}
we can reformulate (\ref{formula_lemmaHyper}) as the statement about partial derivatives of $F$ with respect to the new variable $t$. 
Then, (\ref{formula_lemmaHyper}) can be seen as a claim about the Taylor expansion of $F$ in powers of $t$ and it is equivalent to the claim that 
\beq{
\label{claim_formula}
\frac{2 z^2}{ \Big((1 - t)^2 - 4z\Big)^{3/2}} 
= z^2\sum_{r = 2}^\infty \frac{P_{r - 2} (z) r(r-1)}{(1 - 4z)^{r - 1/2}} t^{r-2} 
= z^2\sum_{n = 0}^\infty\frac{ P_{n}(z)(n+2)(n+1)}{(1 - 4z)^{n + 3/2}} t^{n},
} 
where $P_n(z)$ as defined in (\ref{defi_P}).
In order to derive this formula, we recall that 
\bal{
\frac{1}{(1 - 2x s + s^2)^\alpha} = \sum_{n = 0}^\infty C_n^{(\alpha)}(x) s^n, 
}
where $C_n{(\alpha)}(x)$ are the Gegenbauer polynomials (see  \cite{stein_weiss71}, Section IV.2). Then we use $x = (1 - 4z)^{-1/2}$ and $s = t (1 - 4z)^{-1/2}$, so that 
\bal{
1 - 2xs + s^2 = \frac{(1 - t)^2 - 4z}{1 - 4z},
}
and write
\bal{
\Big((1 - t)^2 - 4z\Big)^{-3/2} = \sum_{n = 0}^\infty C_n^{(3/2)}\Big(\frac{1}{\sqrt{1 - 4z}}\Big) \frac{t^n}{(1 - 4z)^{n/2 + 3/2}}. 
}
Then, we use the fact that 
\bal{
C_n^{(3/2)}\Big(\frac{1}{\sqrt{1 - 4z}}\Big) =
\frac{2^n}{(1 - 4z)^{n/2}} \frac{(3/2)_n}{n!} \tensor*[_2]{F}{_1}\Big(-n/2, (1 - n)/2, -1/2 - n, 1 - 4z\Big),
}
where $(3/2)_n = \Gamma(\frac{3}{2} + n)/\Gamma(\frac{3}{2})$ is the Pochhammer symbol for $3/2$ (see formula 15.9.3 of the Digital Library of Mathematical Functions (\url{https://dlmf.nist.gov/15.9}). It follows that 
\bal{
 F_{uu}(t, z) = 2z^2 \sum_{n = 0}^\infty  \frac{2^n (3/2)_n}{n!} \tensor*[_2]{F}{_1}\Big(-n/2, (1 - n)/2, -1/2 - n, 1 - 4z\Big)\frac{t^n}{(1 - 4z)^{n + 3/2}}.
 }
 It follows that (\ref{claim_formula}) holds with
 \bal{
 P_{n}(z) = \frac{2^{n+1} (3/2)_n}{(n+2)!} \tensor*[_2]{F}{_1}\Big(-n/2, (1 - n)/2, -1/2 - n, 1 - 4z\Big)
 }
for $n \geq 0$, and (\ref{formula_lemmaHyper}) follows by an identity for the hypergeometric functions.

Now we are going to prove (\ref{formula_lemmaHyper_a}). For $r = 1$, it can be verified directly. For $r \geq 2$, (\ref{formula_lemmaHyper_a}) follows from 
(\ref{formula_lemmaHyper}). Indeed, $\beta_{0,r} = 0$ because $P_r(z)$ is a polynomial in $z$ and therefore does not have terms with half-powers of $z$. In order to calculate $\alpha_{0,r}$, we note that formulae (\ref{formula_lemmaHyper}), (\ref{defi_P}) and  (\ref{formula_lemmaHyper_a}) imply that 
  \bal{
 \alpha_{0,r} = \frac{r!}{4^r} P_{r - 2}\Big(\frac{1}{4}\Big) = \Big(\frac{r!}{4^r}\Big) \, \tensor*[_2]{F}{_1}\Big(\frac{3 - r}{2}, 1 -\frac{r}{2}, 2, 1\Big).
 }  
 We use Gauss' identity for the hypergeometric function (see Section 1.3 in \cite{bailey35}): 
 \bal{
 \tensor*[_2]{F}{_1}\Big(a, b, c, 1\Big) = \frac{\Gamma(c) \Gamma(c - a - b)}{\Gamma(c - a) \Gamma(c - b)},
 }
with $c = 2$, $c - a - b = r - \frac{1}{2}$, $c - a = (r+1)/2$, $c - b = r/2 + 1$, which allows us to evaluate:
 \bal{
 \alpha_{0,r} =\frac{ \Gamma(r + 1)}{4^r}\frac{\Gamma(2)\Gamma(r - \frac{1}{2})}{\Gamma(\frac{1}{2} + \frac{r}{2}) \Gamma(1 + \frac{r}{2})} = \frac{\Gamma\big(r - \frac{1}{2}\big)}{2^r \sqrt{\pi}},
 } 
 by using the Legendre duplication identity for the gamma function. This completes the proof.
\end{proof} 
   
%
%

\bigskip
\subsection{Proof of Theorem \ref{theoMainResult0} for $l = 1$}
\label{proofOfMainResult2}

In this section, we will prove a particular case of Theorem \ref{theoMainResult0} in order to illustrate the ``dominant balance'' method for $q$-functional equations.

 From (\ref{exp_derivative}), for $k = 1, 2, \ldots$, we have the following asymptotic expansion around $z = 1/4$:
 \begin{equation}
 \label{exp_Fq}
 \frac{\partial^k}{\partial q^k}  F(q, u, z)\big|_{q = 1, u = 0}=  \frac{\alpha_{k,0}}{(1 - 4z)^{3k/2 - 1/2}} - \frac{\beta_{k,0}}{(1 - 4z)^{3k/2  - 1}} + O\big((1 - 4z)^{-(3k/2  - 3/2)} \big). 
\end{equation}
 In order to determine coefficients $\alpha_{k,0}$ and $\beta_{k,0}$, we change the variables $1 - q = \eps^3$, $1 - 4z = s \eps^2$ (and we set $u = 0$). Then, we 
 define $\wh F(s, \eps) := F(q,u,z)\big|_{u = 0}$ and look for the following expansion of $\wh F(s, \eps)$: 
 \begin{equation}
 \label{exp_dbalance}
 \wh F(s, \eps) = 2 + F_0(s) \eps + F_1(s) \eps^2 + \ldots ,
 \end{equation} 
 where $F_i(s)$ are Laurent series in powers of $s^{1/2}$. More precisely,
 \begin{equation}
 \label{expansions_F0_F1}
 F_0(s) = \sum_{k=0}^\infty \frac{f_{k}^{(0)}}{s^{3k/2 - 1/2}} \text{ and } F_1(s) = \sum_{k=0}^\infty \frac{f_{k}^{(1)}}{s^{3k/2 - 1}}.
 \end{equation}
 For the terms in these expansions we have
 \bal{
 \frac{f_{k}^{(0)}}{s^{3k/2 - 1/2}} \eps &= \frac{f_{k}^{(0)}(1 - q)^{\frac{3k - 1}{3}}}{(1 - 4z)^{\frac{3k - 1}{2}}}(1 - q)^{\frac{1}{3}} = \frac{f_{k}^{(0)}(1 - q)^k}{(1 - 4z)^{\frac{3k - 1}{2}}},
 \\
 \frac{f_{k}^{(1)}}{s^{3k/2 - 1}} \eps^2&=\frac{f_{k}^{(1)}(1 - q)^{\frac{3k - 2}{3}}}{(1 - 4z)^{\frac{3k - 2}{2}}}(1 - q)^{\frac{2}{3}} 
 = \frac{f_{k}^{(1)}(1 - q)^k}{(1 - 4z)^{\frac{3k - 2}{2}}}.
 }
 
 If we differentiate these expressions $k$ times over $q$, set $q = 1$, and compare the result with  (\ref{exp_Fq}), then we find that for $k \geq 0$, we have the relations 
 \begin{equation}
 \label{alpha_beta_0}
 (-1)^k k ! f_{k}^{(0)} = \alpha_{k,0} \text{ and } (-1)^k k! f_{k}^{(1)} = - \beta_{k,0}.
 \end{equation}
 
 In terms of new variables $s$ and $\eps$, the transformation $z \to qz$ can be expressed as follows: 
 \bal{
 s \to s' = \frac{1 - 4qz}{\eps^2} =  \frac{1 - (1 - \eps^3) (1 - s\eps^2)}{\eps^2} = s + \eps(1 - s \eps^2).
 }
 So the functional equation (\ref{eqFun}) can be re-written as follows (after setting $u = 0$): 
 \begin{equation}
 \label{eqFun2}
 \wh F(s, \eps) = \frac{1}{1 - \frac{1}{4}(1 - s\eps^2)\wh F(s + \eps(1 - s\eps^2),\eps)}.
 \end{equation}
 
 We substitute (\ref{exp_dbalance}) in (\ref{eqFun2}) and expand the result in powers of $\eps$. The expansions of $F_i\big(s + \eps(1 - s \eps^2)\big)$ over the small parameter $\eps$ lead to appearance of terms $F'_i(s)$, $F''_i(s)$, and so on. Equating coefficients before $\eps$ and $\eps^2$ on the left and the right hand sides of the expansion of (\ref{eqFun2}) leads to the following differential equations: 
 \begin{align}
 \label{diff_eq_F0}
 F'_0(s) &= 2s - \frac{1}{2} \Big(F_0(s)\Big)^2, \\
 \label{diff_eq_F1}
 F'_1(s) &= F_0(s)\Big(2s - F_1(s) - \frac{1}{2}F'_0(s)\Big) - \frac{1}{2} F''_0(s).
 \end{align} 
 Using (\ref{diff_eq_F0}), equation (\ref{diff_eq_F1}) can be simplified to the following form:
 \begin{equation}
 \label{diff_eq_F1_bis}
  F'_1(s) = F_0(s)\Big(2s - F_1(s)\Big) - 1. 
 \end{equation}
 
 As a remark, equations (\ref{expansions_F0_F1}), (\ref{diff_eq_F0}) and (\ref{diff_eq_F1_bis}) imply the following recursions: $f_0^{(0)} = 2$,  $f_1^{(0)} = -1/2$, $f_0^{(1)} = 2$,  $f_1^{(1)} = -3/2$, and for $n \geq 2$:
 \bal{
 f_n^{(0)} &= \frac{3n - 4}{4} f_{n - 1}^{(0)} - \frac{1}{4} \sum_{i = 1}^{n - 1} f_i^{(0)} f_{n - i}^{(0)},\\
  f_n^{(1)} &=  \frac{3n - 5}{4} f_{n - 1}^{(1)} - \frac{1}{2} \sum_{i = 1}^{n - 1} f_i^{(0)} f_{n - i}^{(1)}.
 }
The first of these recursions is equivalent to the Tak\'acz recursion for the main asymptotic term of the area distribution moments and the second recursion can be used to extract the first order corrections to these moments.
 
 Now we turn to the expansion
 \begin{align}
\frac{\partial^k}{\partial q^k} \frac{\partial}{\partial u} F(q, u, z)\Big|_{q = 1, u = 0}  &= \frac{\alpha_{k,1}}{(1 - 4z)^{3k/2 + 1/2}} - \frac{\beta_{k,1}}{(1 - 4z)^{3k/2 }}
 + O\big((1 - 4z)^{-(3k/2  - 1/2)} \big). \notag
\end{align}
 
We use the same change of variables as before and define 
\bal{
\wh G(s, \eps) := \frac{\partial}{\partial u} F(q, u, z)\Big|_{u = 0}.
} 

From (\ref{solution_q_1}),
\bal{
\frac{\partial}{\partial u} F(q, u, z)\big|_{q = 1, u = 0} = \frac{1}{2}\frac{1}{\sqrt{1 - 4z}}  -\frac{1}{2} = \frac{1}{2}s^{-1/2} \eps^{-1}   -\frac{1}{2} ,
}
and the expansion for $\wh G(s, \eps)$ has the form:
 \begin{equation}
 \label{exp_dbalance2}
 \wh G(s, \eps) = \eps^{-1} (G_0(s) + G_1(s)\eps  + G_2(s) \eps^2 + \ldots) ,
 \end{equation} 
 where $G_i(s)$, $i = 0, 1, 2, \ldots $, can be written as Laurent series in the powers of $s^{1/2}$. In particular, 
  \begin{equation}
 \label{expansions_G0_G1}
 G_0(s) = \sum_{k=0}^\infty \frac{g_{k}^{(0)}}{s^{3k/2 + 1/2}} \text{ and } G_1(s) = \sum_{k=0}^\infty \frac{g_{k}^{(1)}}{s^{3k/2}}.
 \end{equation}
 Then, by an argument similar to the argument before (\ref{alpha_beta_0}), for $k \geq 0$, we have the relations 
 \begin{equation}
  \label{alpha_beta_1}
 (-1)^k k ! g_{k}^{(0)} = \alpha_{k,1} \text{ and } (-1)^k k! g_{k}^{(1)} = -\beta_{k,1}.
 \end{equation}
 
 As the relevant functional equation, we use (\ref{eqFu}), which in the new coordinates and for $u = 0$ can be written as
 \beq{
\label{eqFu1}
\wh G(s, \eps) = \frac{1}{4}(1 - s \eps^2) \Big(\wh F(s, \eps)\Big)^2 \Big(1 +  \wh G(s + \eps(1 - s\eps^2), \eps)\Big).
}
After we substitute (\ref{exp_dbalance}) and (\ref{exp_dbalance2}) in this equation and collect terms before  powers of $\eps$, we obtain several differential equations. The coefficient before $\eps^{-1}$ is identically $0$ and the coefficients before $\eps^0$ and $\eps$ give the following equations:
 \begin{align}
 \label{diff_eq_G0}
 G'_0(s) &=  - \Big(1 +F_0(s) G_0(s)\Big), \\
 \label{diff_eq_G1}
 G'_1(s) &= -\frac{1}{4} F_0^2(s) G_0(s) + \Big(s - F_1(s)\Big)  G_0(s) - F_0(s)\Big(1 + G_1(s) + G'_0(s)\Big) - \frac{1}{2} G''_0(s).
 \end{align} 
 The second equation can be slightly simplified by using the first one:
 \beq{
  \label{diff_eq_G1_bis}
 G'_1(s) = -\frac{3}{4} F_0(s)\Big( 1+  G'_0(s)\Big) + \Big(s - F_1(s)\Big)  G_0(s) - F_0(s) G_1(s) - \frac{1}{2} G''_0(s).
 }
 These equations can be easily converted to recursions for $g_{k}^{(0)}$ and $g_{k}^{(1)}.$ However, we will use a different approach in order to relate $g_{k}^{(0)}$ and $g_{k}^{(1)}$ to $f_{k}^{(0)}$ and $f_{k}^{(1)}$, respectively.  

 \begin{lemma}
 \label{lemma_connection}
 $G_0(s) = -\frac{1}{2} F'_0(s)$ and $G_1(s) = \frac{1}{2} -\frac{1}{2} F'_1(s)$.
 \end{lemma}

 \begin{proof}
 If we try $-\frac{1}{2} F'_0(s)$ instead of $G_0(s)$ in (\ref{diff_eq_G0}), we get 
 \bal{
 -\frac{1}{2} F''_0(s) &= - \Big(1 - F_0(s)\frac{1}{2} F'_0(s)\Big) = -1 + \frac{1}{4} \Big(F_0(s)^2\Big)', \text{ or }
 \\
 F''_0(s) &= 2 - \frac{1}{2} \Big(F_0(s)^2\Big)',
 }
 which is an immediate consequence of  (\ref{diff_eq_F0}). It follows  that  $-\frac{1}{2}F'_0(s)$ satisfies the same differential equation as $G_0(s)$. In addition, we can check that  the leading term in the expansion (\ref{expansions_F0_F1}) for $F_0(s)$ is $-2 s^{1/2}$ and the leading term in the expansion (\ref{expansions_G0_G1}) for $G_0(s)$ is $\frac{1}{2}s^{-1/2}$. 
 
 Indeed, the leading term for $F_0(z)$ is $f_0^{(0)}s^{1/2} = \alpha_{0, 0} s^{1/2}$ by  (\ref{expansions_F0_F1}) and (\ref{alpha_beta_0}), and by defining expansion (\ref{exp_derivative}), $\alpha_{0, 0}$ is the coefficient before $\sqrt{1 - 4z}$ in the expansion of $D(q, u, z)|_{q = 1, u = 0}$ in powers of $1 - 4z$. From (\ref{equ_Fz}), we can read it off as $\alpha_{0, 0} = - 2$. 
 Then, the leading term for $G_0(s)$ is $g_0^{(0)} s^{-1/2} = \alpha_{0, 1} s^{-1/2} = \frac{1}{2} s^{-1/2}$ by Proposition  \ref{propoAlpha}.
 
  This imply that $G_0(s)$ and $-\frac{1}{2} F'_0(s)$ have the same leading terms. Together with the fact that these functions satisfy the same differential equation (\ref{diff_eq_G0}), this implies that $G_0(s) = - \frac{1}{2} F'_0(s)$.

 The second claim is proved similarly by substituting $ -\frac{1}{2} F'_0(s)$ and $\frac{1}{2} -\frac{1}{2} F'_1(s)$ instead of $G_0(s)$ and $G_1(s)$in (\ref{diff_eq_G1_bis}) and using equations  (\ref{diff_eq_F0}) and (\ref{diff_eq_F1_bis}) to eliminate the derivatives of $F_0(s)$ and $F_1(s)$. This shows that $\frac{1}{2} -\frac{1}{2} F'_1(s)$ satisfies the same differential equation as $G_1(s)$. Then one can check the equality of the leading terms in the corresponding expansions.  
 
 Indeed, $F_1(s) = f_0^{(1)} s + \ldots = -\beta_{0, 0} s + \ldots$, and $-\beta_{0, 0}$ is the coefficient before $(1 - 4z)$ in the expansion of $D(q, u, z)|_{q = 1, u = 0}$ in powers of $1 - 4z$. From (\ref{equ_Fz}), we read it off as $-\beta_{0, 0} = 2$. Then, $G_1(s) = g_0^{(1)} + g_1^{(1)} s^{-3/2} + \ldots$ and $g_0^{(1)} = - \beta_{0, 1} = - \frac{1}{2}$ from Proposition \ref{propoAlpha}. Therefore, $G_1(s)$ and $\frac{1}{2} - \frac{1}{2} F'_1(s)$ have the same leading term and satisfy the same differential equation, so we can conclude that $G_1(s) = \frac{1}{2} - \frac{1}{2} F'_1(s)$. 
  \end{proof}

 Equations $G_0(s) = - \frac{1}{2}F'_0(s)$ and $G_1(s) = \frac{1}{2} -\frac{1}{2} F'_1(s)$, and expansions (\ref{expansions_F0_F1}) and (\ref{expansions_G0_G1}) imply that for all $k \geq 1$, 
 \bal{
 g_k^{(0)} = \frac{3k - 1}{4} f_k^{(0)}, \\
 g_k^{(1)} = \frac{3k - 2}{4} f_k^{(1)}.
 }
 Then, in terms of coefficients $\alpha$ and $\beta$, the second of these equations can be re-written as
 \begin{equation}
 \label{main_identity4bis}
 \beta_{k,1} = \frac{3k - 2}{4} \beta_{k,0},
 \end{equation}
 where we used  (\ref{alpha_beta_0}) and  (\ref{alpha_beta_1}). This is exactly the identity (\ref{main_identity4}) that we needed to prove Theorem \ref{theoMainResult0} for the case $l = 1$.

%
%

\subsection{Proof of Theorem \ref{theoMainResult0}}
\label{proofOfAsymptoticIndependence}
We have the expansion
 \begin{align}
 \label{expansion_basic}
\frac{\partial^k}{\partial q^k} \frac{\partial^l}{\partial u^l} F(q, u, z)\Big|_{q = 1, u = 0}  &= \frac{\alpha_{k,l}}{(1 - 4z)^{3k/2 + l - 1/2}} - \frac{\beta_{k,l}}{(1 - 4z)^{3k/2 +l - 1 }}\\
& + O\big((1 - 4z)^{-(3k/2  + l - 3/2)} \big). \notag
\end{align}
We use the following scaling: $1 - q = \eps^3$, $1 - 4z = s \eps^2$, $u = t \eps^2$ and define 
\bal{
\wh F(s, t, \eps) &= F(q, u, z).
}
Note that  the transformation $z \to qz$ can be expressed as follows: 
 \bal{
 s \to s' = \frac{1 - 4qz}{\eps^2} =  \frac{1 - (1 - \eps^3) (1 - s\eps^2)}{\eps^2} = s + \eps(1 - s \eps^2).
 }
 Then the functional equation (\ref{eqFun}) can be written as
 \beq{
\label{eqFu_bis}
\wh F(s, t, \eps) = \frac{1}{1 - \frac{1}{4} t\eps^2(1 - s \eps^2) - \frac{1}{4} (1 - s \eps^2) \wh F\big(s + \eps(1 - s\eps^2), t, \eps\big) }.
}
We look for solution in the form:
\beq{
\label{approximation1}
\wh F(s, t, \eps) = 2 + (2s - t/2) \eps^2 + \eps\Big(F^0(s, t)+ F^1(s,t) \eps + F^2(s, t) \eps^2 + \ldots \Big).
}
(The term $(2s - t/2)\eps^2$ is for convenience in some further calculations.) Before analyzing (\ref{eqFu_bis}) and (\ref{approximation1}), let us explain what is the ultimate goal. 

Since the change of variables $1 - q = \eps^3$, $1 - 4z = s \eps^2$, $u = t \eps^2$ implies that  
\bal{
\frac{\partial^l f(q, u , z)}{\partial t^l}  = \eps^{2l} \frac{\partial^l  f(q, u, z)}{\partial u^l},
} 
we have for $l \geq 1$, 
\begin{equation}
 \label{approximation0}
\frac{\partial^l}{\partial u^l} F(q, u, z)\Big|_{u = 0} = \eps^{1-2l}\Big(F_{l, 0}(s) + F_{l, 1}(s) \eps +  F_{l, 2}(s) \eps^2 + \ldots\Big),
\end{equation}
 where
\beq{
\label{relation_approximations}
 F_{l, i}(s) = \frac{\partial^l}{\partial t^l} \Big((2s - t/2)\delta_{i,1}  + F^i(s, t)\Big)\Big|_{t = 0},
}
with $\delta_{i,1} = 1$ if $i = 1$ and zero otherwise.

The coefficients in this expansion have Laurent series
\beq{
\label{expansion_F_li}
 F_{l,i}(s) = \sum_{k = 0}^\infty \frac{f_{k, l, i}}{s^{3k/2 + l -  1/2 - i/2}},
}
Consider a specific term in the expansion (\ref{approximation0}):
\bal{
F_{l, i} \eps^{i + 1 - 2l} &= \sum_{k = 0}^\infty \frac{f_{k, l, i}}{s^{3k/2 + l - 1/2 - i/2}} \eps^{i + 1 - 2l} = \sum_{k = 0}^\infty \frac{f_{k, l, i}}{\big((1-4z)/\eps^2\big)^{3k/2 + l - 1/2 - i/2}} \eps^{i + 1 - 2l} \\
&= \sum_{k = 0}^\infty \frac{f_{k, l, i}}{(1-4z)^{3k/2 + l - 1/2 - i/2}} \eps^{3k} = \sum_{k = 0}^\infty \frac{f_{k, l, i}}{(1-4z)^{3k/2 + l - 1/2 - i/2}} (1 - q)^{k}.
}
 Differentiate this series $k$-times with respect to $q$, set $q = 0$, and compare with (\ref{expansion_basic}). We find that 
\begin{align}
\notag
\alpha_{k,l} &= (-1)^k k! f_{k, l, 0}, \\
\label{alpha_beta_f}
-\beta_{k,l} &= (-1)^k k! f_{k, l, 1},
\end{align}
valid for $k, l \geq 1$. We are looking for relations between $f_{k, l, 1}$ for different $l$, since this would give us  the desired relation between coefficients $\beta_{k,l}$ for different $l$. More generally, we will look for relations between $F_{l,1}(1)$, or, in other words $(\partial^l t) F^{1}(s, t)$ for different $l$. We will also need relations between $(\partial^l t) F^{0}(s, t)$.

We substitute expansion (\ref{approximation1}) in equation (\ref{eqFu_bis}) and expand the result in powers of $\eps$. Then, we find that for powers $0$ and $1$ the equation is satisfied trivially, and for powers $2$ and $3$, we get the following differential equations satisfied by $F^0(s, t)$ and $F^1(s,t)$:
\begin{align}
\label{diff_eq_F0bis}
\partial_s F^0 &= 2s - t - \frac{1}{2} \Big(F^0\Big)^2,\\
\label{diff_eq_F1bis}
\partial_s F^1 &= -F^0 F^1 - \frac{1}{2}F^0\partial_s F^0 - \frac{1}{2}\partial_{ss} F^0.
\end{align}

%
\begin{propo}
\label{propo_derivatives_F0}
For all $k \geq 1$, we have 
\beq{
\frac{\partial^k}{\partial t^k} F^0(s,t) = \Big(-\frac{1}{2}\Big)^k  \frac{\partial^k}{\partial s^k} F^0(s,t).
}
\end{propo}

\begin{proof}
Let $F_t^0(s,t) = \frac{\partial}{\partial t} F^0(s,t)$ and $F_s^0(s,t) = \frac{\partial}{\partial s} F^0(s,t)$. (To make the notation less cumbersome, we omit the superscript $0$ for function $F^0(s, t)$, $F_t^0(s, t)$, and $F_s^1(s,t)$ in the following.) 
If we differentiate equation (\ref{diff_eq_F0bis}) either by $t$ or by $s$, we get differential equations for these functions:
\begin{align}
\label{equ_first_derivatives1}
\partial_s F_{t} &= - 1 - 2F_t F, \text{ and }
\\
\label{equ_first_derivatives2}
 \partial_s F_{s} &= 2 - 2F_s F, 
\end{align}
respectively. The functions $F_t$ and $F_s$ are determined by these equations and initial conditions.  
If we try $-\frac{1}{2}F_s$ as a possible solution of (\ref{equ_first_derivatives1}), i.e., if we substitute $ F_t = -\frac{1}{2} F_s$ in (\ref{equ_first_derivatives1}), then it becomes:
\bal{
-\frac{1}{2}\partial_s F_{s} = - 1 +  F_s F,
} 
which holds because it is the consequence of (\ref{equ_first_derivatives2}). It follows that $-\frac{1}{2} F_s$ satisfies the same differential equation (\ref{equ_first_derivatives1}) as $F_t$.

In addition,  we have formula (\ref{solution_q_1}) for $F(q, u, z)|_{q = 1}$. If we do a change of variables $(z, u, q) \to (s, t, \eps)$ and expand the result in powers of $\eps$, then we get:
\beq{
\label{singular_solution}
F_{sing}(s, t, \eps) = 2 - \sqrt{4 s - 2 t} \eps + (2s - t/2)\eps^2 + \ldots,
}
which implies that the equation $F_t^0= -\frac{1}{2} F_s^0$ holds at the surface $q = 1$ near the singular point for  function $F^0(s, t) =\sqrt{4s - 2t}$. So, this equality holds in general. The validity of this equation shows that the claim of the proposition holds for $k = 1$. 

We proceed by induction over $k$. Differentiating (\ref{equ_first_derivatives1}) further $k-1$ times over $t$ or  $s$, respectively, we obtain the expressions
\beq{
\label{equ_2}
\partial_s \partial_t^k F = \sum_{\substack{0 \leq i \leq j \leq k\\ i + j = k}}^k a_{i, j}(\partial_t^i F) (\partial_t^{j} F)
\text{ and } \partial_s \partial_s^k F = \sum_{\substack{0 \leq i \leq j \leq k\\ i + j = k}} a_{i, j}(\partial_s^i F) (\partial_s^{j} F), 
} 
where, by convention, $\partial_s^0 F = \partial_t^0 F := F$, and $a_{i, j}$ are some coefficients, with $a_{0,k} = - 2$. 
For the inductive step, assume that $\partial_t^{i} F = \big(-\frac{1}{2}\big)^i \partial_s^{i} F$ for all $i < k$. Then, the first of the above equations becomes:
\beq{
\label{equ_3}
\partial_s \partial_t^k F = -2 F (\partial_t^k F ) + \sum_{\substack{0 \leq i \leq j \leq k-1\\ i + j = k}} a_{i, j}\Big(-\frac{1}{2}\Big)^k (\partial_s^i F) (\partial_s^{j} F).
}
We consider this equation as an equation for function $\partial_t^k F$. It is clear that $\big(-\frac{1}{2}\big)^k \partial_s^{k} F$ also satisfies this equation since upon  substitution of  $\partial_t^k F:=\big(-\frac{1}{2}\big)^k \partial_s^{k} F$, the equation (\ref{equ_3}) becomes a consequence of the second equation in (\ref{equ_2}). Since $\partial_t^k F:=\big(-\frac{1}{2}\big)^k \partial_s^{k} F$ holds also for the singular solution (\ref{singular_solution}), we conclude that it holds in general. 
\end{proof}

%
\begin{propo}
\label{propo_derivatives_F1}
For all $k \geq 1$, we have
\bal{
\frac{\partial^k}{\partial t^k} F^1(s,t) = \Big(-\frac{1}{2}\Big)^k  \frac{\partial^k}{\partial s^k} F^1(s,t).
}
\end{propo}
\begin{proof}
Here we start with equation (\ref{diff_eq_F1bis}). Differentiating it over $t$ and $s$, respectively, we find
\beq{
\label{equ_4}
\partial_s F_t^1 = -F_t^0 F^1 - F^0 F^1_t - \frac{1}{2} F^0_t F^0_s - \frac{1}{2} F^0 F^0_{t,s}- \frac{1}{2} F^0_{t,ss}
} 
and
\beq{
\label{equ_5}
\partial_s F_s^1 = -F_s^0 F^1 - F^0 F^1_s - \frac{1}{2} F^0_s F^0_s - \frac{1}{2} F^0 F^0_{ss} - \frac{1}{2} F^0_{sss}.
}

From Proposition \ref{propo_derivatives_F0}, we have $F_t^0 = -\frac{1}{2} F_s^0$. We plug this identity in equation (\ref{equ_4}) and find:
\beq{
\label{equ_6}
-2 \partial_s F_t^1 = -F_s^0 F^1 + 2 F^0 F^1_t - \frac{1}{2} F^0_s F^0_s - \frac{1}{2} F^0 F^0_{ss}- \frac{1}{2} F^0_{sss}.
}

We are looking at this equation as a differential equation for $F_t^1$ and try $F_t^1 := - \frac{1}{2}F_s^1$ as a possible solution. We plug this function in equation (\ref{equ_6}), and see that it it is a valid solution because of equation (\ref{equ_5}). Referring to (\ref{singular_solution}), we find that this equality holds also for $F_{sing}^1 = 0$. We conclude that it holds in general. 

For the case $k \geq 2$ we differentiate further. We differentiate equations of type (\ref{equ_5}) over $s$ and equations of type (\ref{equ_6}) over $t$. While differentiating over $t$, we use the identities $\partial_t^k F^0 = (-1/2)^k \partial_s^k F^0$, and the identities $\partial_t^j F^1 = (-1/2)^j \partial_s^j F^1$, which by inductive assumption are known for $j < k$. As a result, at every step, we obtain equations of the form:
\beq{
\label{equ_7}
\partial_s^{k+1} F^1 = - F^0 (\partial^{k}_s F^1) - \sum_{i=0}^{k - 1} a_{i, j} (\partial^i_s F^1) (\partial^{k - i}_s F^0) - \sum_{\substack{0 \leq i \leq j \leq k\\ i + j = k}} b_{i,j} (\partial^i_s F^0) (\partial^{j}_s F^0)
- \frac{1}{2} \partial^{k + 2}_s F^0, 
}
and
\begin{align}
\label{equ_8}
&(-2)^k\partial_s(\partial_t^{k} F^1) =\\
&= -(-2)^k F^0 (\partial^{k}_t F^1) - \sum_{i=0}^{k - 1} a_{i, j} (\partial^i_s F^1) (\partial^{k - i}_s F^0) - \sum_{\substack{0 \leq i \leq j \leq k\\ i + j = k}} b_{i,j} (\partial^i_s F^0) (\partial^{j}_s F^0)
- \frac{1}{2} \partial^{k + 2}_s F^0, \notag
\end{align}
and we observe that $\partial_t^{k} F^1$ defined as $(-1/2)^k\partial_s^{k} F^1$ satisfies the second equation because it reduces this equation to the first one. (In the next step, we differentiate (\ref{equ_7}) over $s$ and (\ref{equ_8}) over $t$, multiply the result of the latter by $-2$ and proceed to the next inductive step.)
\end{proof}

We use relations (\ref{expansion_F_li}), (\ref{alpha_beta_f})  and (\ref{relation_approximations}) in order to rewrite Proposition \ref{propo_derivatives_F1} in the following form:
\beq{
\sum_{k = 0}^\infty \frac{\beta_{k, l}/k!}{s^{3k/2 + l -  1}} = \Big(-\frac{1}{2}\Big)^l  \frac{\partial^l}{\partial s^l} \bigg(\sum_{k = 0}^\infty \frac{\beta_{k, 0}/k!}{s^{3k/2 -  1}}\bigg).
}
Here we used $i = 1$ and $l \geq 2$. (The case $l = 1$ was handled in the previous section.) In terms of coefficients $\beta_{k, l}$ this means 
\beq{
\label{main_identity2a_bis}
\beta_{k,l} = \Big(\frac{1}{2}\Big)^l \Big(\frac{3k}{2} - 1\Big) \frac{3k}{2}\times \ldots \times\Big(\frac{3k}{2} + l - 2\Big) \beta_{k, 0}.
}
This is the identity (\ref{main_identity2a}) that we needed to prove the general case of Theorem \ref{theoMainResult0}.


\begin{thebibliography}{31}

\bibitem[\protect\citeauthoryear{Bailey}{1935}]{bailey35}
\begin{bbook}[author]
\bauthor{\bsnm{Bailey},~\bfnm{W.~N.}\binits{W.~N.}}
(\byear{1935}).
\btitle{Generalized Hypegeometric Series}.
\bpublisher{Cambridge University Press}.
\end{bbook}
\endbibitem

\bibitem[\protect\citeauthoryear{Berndt and Yee}{2003}]{berndt_yee2003}
\begin{barticle}[author]
\bauthor{\bsnm{Berndt},~\bfnm{B.~C.}\binits{B.~C.}} \AND
  \bauthor{\bsnm{Yee},~\bfnm{A.~J.}\binits{A.~J.}}
(\byear{2003}).
\btitle{On the generalized \mbox{R}ogers-\mbox{R}amanujan continued fraction}.
\bjournal{Ramanujan J.}
\bvolume{7}
\bpages{321--331}.
\end{barticle}
\endbibitem

\bibitem[\protect\citeauthoryear{Blanco and
  Petersen}{2014}]{blanco_petersen2014}
\begin{barticle}[author]
\bauthor{\bsnm{Blanco},~\bfnm{Saul~A.}\binits{S.~A.}} \AND
  \bauthor{\bsnm{Petersen},~\bfnm{T.~Kyle}\binits{T.~K.}}
(\byear{2014}).
\btitle{Counting \mbox{D}yck paths by area and rank}.
\bjournal{Ann. Comb.}
\bvolume{18}
\bpages{171--197}.
\end{barticle}
\endbibitem

\bibitem[\protect\citeauthoryear{Chassaing, Marckert and
  Yor}{2000}]{chassaing_marckert_yor2000}
\begin{binproceedings}[author]
\bauthor{\bsnm{Chassaing},~\bfnm{Philippe}\binits{P.}},
  \bauthor{\bsnm{Marckert},~\bfnm{Jean-Fran\c{c}ois}\binits{J.-F.}} \AND
  \bauthor{\bsnm{Yor},~\bfnm{Marc}\binits{M.}}
(\byear{2000}).
\btitle{The height and width of simple trees}.
In \bbooktitle{Mathematics and Computer Science}
(\beditor{\bfnm{Dani\'{e}le}\binits{D.}~\bsnm{Gardy}} \AND
  \beditor{\bfnm{Abdelkader}\binits{A.}~\bsnm{Mokkadem}}, eds.)
\bpages{17--30}.
\bpublisher{Birkh{\"a}user Basel}.
\end{binproceedings}
\endbibitem

\bibitem[\protect\citeauthoryear{Cs\"org\H{o}, Shi and
  Yor}{1999}]{csorgo_shi_yor99}
\begin{barticle}[author]
\bauthor{\bsnm{Cs\"org\H{o}},~\bfnm{Mikl\'os}\binits{M.}},
  \bauthor{\bsnm{Shi},~\bfnm{Zhan}\binits{Z.}} \AND
  \bauthor{\bsnm{Yor},~\bfnm{Marc}\binits{M.}}
(\byear{1999}).
\btitle{Some asymptotic properties of the local time of the uniform empirical
  process}.
\bjournal{Bernoulli}
\bvolume{5}
\bpages{1035--1058}.
\end{barticle}
\endbibitem

\bibitem[\protect\citeauthoryear{Drmota}{2009}]{drmota2009}
\begin{bbook}[author]
\bauthor{\bsnm{Drmota},~\bfnm{Michael}\binits{M.}}
(\byear{2009}).
\btitle{Random Trees}.
\bpublisher{Springer Verlag/Wien}.
\end{bbook}
\endbibitem

\bibitem[\protect\citeauthoryear{Fill and Janson}{2009}]{fill_janson2009}
\begin{barticle}[author]
\bauthor{\bsnm{Fill},~\bfnm{James}\binits{J.}} \AND
  \bauthor{\bsnm{Janson},~\bfnm{Svante}\binits{S.}}
(\byear{2009}).
\btitle{Precise logarithmic asymptotics for the right tails of some limit
  random variables for random trees}.
\bjournal{Ann. Comb.}
\bvolume{12}
\bpages{403--416}.
\end{barticle}
\endbibitem

\bibitem[\protect\citeauthoryear{Flajolet and
  Louchard}{2001}]{flajolet_louchard2001}
\begin{barticle}[author]
\bauthor{\bsnm{Flajolet},~\bfnm{Philippe}\binits{P.}} \AND
  \bauthor{\bsnm{Louchard},~\bfnm{Guy}\binits{G.}}
(\byear{2001}).
\btitle{Analytic variations on the \mbox{A}iry distribution}.
\bjournal{Algorithmica}
\bvolume{31}
\bpages{361--377}.
\end{barticle}
\endbibitem

\bibitem[\protect\citeauthoryear{Flajolet and
  Sedgewick}{2009}]{flajolet_sedgewick2009}
\begin{bbook}[author]
\bauthor{\bsnm{Flajolet},~\bfnm{Philippe}\binits{P.}} \AND
  \bauthor{\bsnm{Sedgewick},~\bfnm{Robert}\binits{R.}}
(\byear{2009}).
\btitle{Analytic Combinatorics}.
\bpublisher{Cambridge University Press}.
\end{bbook}
\endbibitem

\bibitem[\protect\citeauthoryear{Janson}{2001}]{janson2001}
\begin{barticle}[author]
\bauthor{\bsnm{Janson},~\bfnm{Svante}\binits{S.}}
(\byear{2001}).
\btitle{Moment convergence in conditional limit theorems}.
\bjournal{J. Appl. Probab.}
\bvolume{38}
\bpages{421--438}.
\end{barticle}
\endbibitem

\bibitem[\protect\citeauthoryear{Janson}{2007}]{janson2007}
\begin{barticle}[author]
\bauthor{\bsnm{Janson},~\bfnm{Svante}\binits{S.}}
(\byear{2007}).
\btitle{Brownian excursion area, \mbox{W}right\mbox{'}s constants in graph
  enumeration, and other \mbox{B}rownian areas}.
\bjournal{Probab. Surv.}
\bvolume{4}
\bpages{80--145}.
\end{barticle}
\endbibitem

\bibitem[\protect\citeauthoryear{Janson}{2008}]{janson2008}
\begin{barticle}[author]
\bauthor{\bsnm{Janson},~\bfnm{Svante}\binits{S.}}
(\byear{2008}).
\btitle{On the asymptotic joint distribution of height and width in random
  trees}.
\bjournal{Stud. Sci. Math. Hung.}
\bvolume{45}
\bpages{451--467}.
\end{barticle}
\endbibitem

\bibitem[\protect\citeauthoryear{Kolchin}{1984}]{kolchin84}
\begin{bbook}[author]
\bauthor{\bsnm{Kolchin},~\bfnm{V.~F.}\binits{V.~F.}}
(\byear{1984}).
\btitle{Random Mappings},
\bedition{\mbox{R}ussian} ed.
\bseries{(English translation: Optimization Software, New York, 1986)}.
\bpublisher{Nauka, Moscow}.
\end{bbook}
\endbibitem

\bibitem[\protect\citeauthoryear{Labarbe and
  Marckert}{2007}]{labarbe_marckert2007}
\begin{barticle}[author]
\bauthor{\bsnm{Labarbe},~\bfnm{Jean-Maxime}\binits{J.-M.}} \AND
  \bauthor{\bsnm{Marckert},~\bfnm{Jean-Fran\c{c}ois}\binits{J.-F.}}
(\byear{2007}).
\btitle{Asymptotics of \mbox{B}ernoulli random walks, bridges, excursions and
  meanders with a given number of peaks}.
\bjournal{Electron. J. Probab.}
\bvolume{12}
\bpages{229--261}.
\end{barticle}
\endbibitem

\bibitem[\protect\citeauthoryear{Louchard}{1984a}]{louchard1984}
\begin{barticle}[author]
\bauthor{\bsnm{Louchard},~\bfnm{G.}\binits{G.}}
(\byear{1984}a).
\btitle{Kac\mbox{'}s formula, \mbox{L}evy\mbox{'}s local time and
  \mbox{B}rownian excursion}.
\bjournal{J. Appl. Probab.}
\bvolume{21}
\bpages{479--499}.
\end{barticle}
\endbibitem

\bibitem[\protect\citeauthoryear{Louchard}{1984b}]{louchard1984b}
\begin{barticle}[author]
\bauthor{\bsnm{Louchard},~\bfnm{G.}\binits{G.}}
(\byear{1984}b).
\btitle{The \mbox{B}rownian excursion area: a numerical analysis}.
\bjournal{Comp. \mbox{\&} Maths. with Appls.}
\bvolume{10}
\bpages{413--417}.
\end{barticle}
\endbibitem

\bibitem[\protect\citeauthoryear{Narayana}{1959}]{narayana59}
\begin{barticle}[author]
\bauthor{\bsnm{Narayana},~\bfnm{T.~V.}\binits{T.~V.}}
(\byear{1959}).
\btitle{A partial order and its applications to probability theory}.
\bjournal{Sanky\mbox{\=a}}
\bvolume{21}
\bpages{91--98}.
\end{barticle}
\endbibitem

\bibitem[\protect\citeauthoryear{{Ngu\mbox{\~{\^e}}n
  Th\mbox{\'{\^e}}}}{2004}]{the2004}
\begin{barticle}[author]
\bauthor{\bsnm{{Ngu\mbox{\~{\^e}}n
  Th\mbox{\'{\^e}}}},~\bfnm{Michel}\binits{M.}}
(\byear{2004}).
\btitle{Area and inertial moment of \mbox{D}yck paths}.
\bjournal{Comb. Probab. Comput.}
\bvolume{13}
\bpages{697--716}.
\end{barticle}
\endbibitem

\bibitem[\protect\citeauthoryear{Petersen}{2015}]{petersen2015}
\begin{bbook}[author]
\bauthor{\bsnm{Petersen},~\bfnm{T.~Kyle}\binits{T.~K.}}
(\byear{2015}).
\btitle{Eulerian Numbers}.
\bpublisher{Birkhauser}.
\end{bbook}
\endbibitem

\bibitem[\protect\citeauthoryear{Richard}{2002}]{richard2002}
\begin{barticle}[author]
\bauthor{\bsnm{Richard},~\bfnm{Christoph}\binits{C.}}
(\byear{2002}).
\btitle{Scaling behaviour of two-dimensional polygon models}.
\bjournal{J. Stat. Phys.}
\bvolume{108}
\bpages{459--493}.
\end{barticle}
\endbibitem

\bibitem[\protect\citeauthoryear{Richard}{2009}]{richard2009}
\begin{barticle}[author]
\bauthor{\bsnm{Richard},~\bfnm{Christoph}\binits{C.}}
(\byear{2009}).
\btitle{On $q$-functional equations and excursion moments}.
\bjournal{Discrete Math.}
\bvolume{309}
\bpages{207--230}.
\end{barticle}
\endbibitem

\bibitem[\protect\citeauthoryear{Richard and
  Guttmann}{2001}]{richard_guttmann2001}
\begin{barticle}[author]
\bauthor{\bsnm{Richard},~\bfnm{Christoph}\binits{C.}} \AND
  \bauthor{\bsnm{Guttmann},~\bfnm{Anthony~J.}\binits{A.~J.}}
(\byear{2001}).
\btitle{$q$-linear approximants: Scaling functions for polygon models}.
\bjournal{J. Phys. A Math.}
\bvolume{34}
\bpages{4783--4796}.
\end{barticle}
\endbibitem

\bibitem[\protect\citeauthoryear{Richard, Guttmann and
  Jensen}{2001}]{richard_guttmann_jensen2001}
\begin{barticle}[author]
\bauthor{\bsnm{Richard},~\bfnm{Christoph}\binits{C.}},
  \bauthor{\bsnm{Guttmann},~\bfnm{Anthony~J.}\binits{A.~J.}} \AND
  \bauthor{\bsnm{Jensen},~\bfnm{Iwan}\binits{I.}}
(\byear{2001}).
\btitle{Scaling function and universal amplitude combinations for self-avoiding
  polygons}.
\bjournal{J. Phys. A Math.}
\bvolume{34}
\bpages{L495--501}.
\end{barticle}
\endbibitem

\bibitem[\protect\citeauthoryear{Richard, Jensen and
  Guttmann}{2008}]{richard_jensen_guttmann2008}
\begin{barticle}[author]
\bauthor{\bsnm{Richard},~\bfnm{Christoph}\binits{C.}},
  \bauthor{\bsnm{Jensen},~\bfnm{Iwan}\binits{I.}} \AND
  \bauthor{\bsnm{Guttmann},~\bfnm{Anthony~J.}\binits{A.~J.}}
(\byear{2008}).
\btitle{Area distribution and scaling function for punctured polygons}.
\bjournal{Electron. J. Comb.}
\bvolume{15}
\bpages{R53}.
\end{barticle}
\endbibitem

\bibitem[\protect\citeauthoryear{Schwerdtfeger, Richard and
  Thatte}{2010}]{srt2010}
\begin{barticle}[author]
\bauthor{\bsnm{Schwerdtfeger},~\bfnm{Uwe}\binits{U.}},
  \bauthor{\bsnm{Richard},~\bfnm{Christoph}\binits{C.}} \AND
  \bauthor{\bsnm{Thatte},~\bfnm{Bhalchandra}\binits{B.}}
(\byear{2010}).
\btitle{Area limit laws for symmetry classes of staircase polygons}.
\bjournal{Comb. Probab. Comput.}
\bvolume{19}
\bpages{441--461}.
\end{barticle}
\endbibitem

\bibitem[\protect\citeauthoryear{Stanley}{1999}]{stanley99}
\begin{bbook}[author]
\bauthor{\bsnm{Stanley},~\bfnm{Richard~P.}\binits{R.~P.}}
(\byear{1999}).
\btitle{Enumerative Combinatorics, volume 2}.
\bpublisher{Cambridge University Press}.
\end{bbook}
\endbibitem

\bibitem[\protect\citeauthoryear{Stanley}{2015}]{stanley2015}
\begin{bbook}[author]
\bauthor{\bsnm{Stanley},~\bfnm{Richard~P.}\binits{R.~P.}}
(\byear{2015}).
\btitle{Catalan Numbers}.
\bpublisher{Cambridge University Press}.
\end{bbook}
\endbibitem

\bibitem[\protect\citeauthoryear{Stein and Weiss}{1971}]{stein_weiss71}
\begin{bbook}[author]
\bauthor{\bsnm{Stein},~\bfnm{Elias}\binits{E.}} \AND
  \bauthor{\bsnm{Weiss},~\bfnm{Guido}\binits{G.}}
(\byear{1971}).
\btitle{Introduction to Fourier Analysis on Euclidean Spaces}.
\bpublisher{Princeton University Press: Princeton, New Jercey}.
\end{bbook}
\endbibitem

\bibitem[\protect\citeauthoryear{Tak\'{a}cs}{1991}]{takacs91}
\begin{barticle}[author]
\bauthor{\bsnm{Tak\'{a}cs},~\bfnm{Lajos}\binits{L.}}
(\byear{1991}).
\btitle{A \mbox{B}ernoulli excursion and its various applications}.
\bjournal{Adv. Appl. Probab.}
\bvolume{23}
\bpages{557--585}.
\end{barticle}
\endbibitem

\bibitem[\protect\citeauthoryear{Tak\'{a}cs}{1992}]{takacs1992}
\begin{barticle}[author]
\bauthor{\bsnm{Tak\'{a}cs},~\bfnm{Lajos}\binits{L.}}
(\byear{1992}).
\btitle{Random walk processes and their applications in order statistics}.
\bjournal{Ann. Appl. Probab.}
\bvolume{2}
\bpages{435--459}.
\end{barticle}
\endbibitem

\bibitem[\protect\citeauthoryear{Tak\'{a}cs}{1994}]{takacs94}
\begin{barticle}[author]
\bauthor{\bsnm{Tak\'{a}cs},~\bfnm{Lajos}\binits{L.}}
(\byear{1994}).
\btitle{On the total height of random rooted binary trees}.
\bjournal{J. Comb. Theory Ser. B}
\bvolume{61}
\bpages{155--166}.
\end{barticle}
\endbibitem

\end{thebibliography}
\end{document}